\documentclass[12pt]{amsart}
 \usepackage[mathscr]{eucal}
 \usepackage[all]{xypic}
\author{Charles Cadman}
\title{Gromov-Witten Invariants of \boldmath $\bP^2$-Stacks}
\thanks{to appear in {\sl Compositio Mathematica}.  This version differs from the published version.}

\numberwithin{equation}{subsection}
\theoremstyle{plain}
\newtheorem{thm}[equation]{Theorem}
\newtheorem{lemma}[equation]{Lemma}
\newtheorem{prop}[equation]{Proposition}
\newtheorem{cor}[equation]{Corollary}
\theoremstyle{definition}
\newtheorem{definition}[equation]{Definition}
\theoremstyle{remark}
\newtheorem{remark}[equation]{Remark}
\newtheorem{algorithm}[equation]{Algorithm}

\newcommand{\Der}{\mathrm{Der}}
\newcommand{\Sym}{\mathrm{Sym}}
\newcommand{\Tor}{\mathrm{Tor}}
\newcommand{\spec}{\mathrm{Spec}\;}
\newcommand{\eC}{\mathscr{C}}
\newcommand{\eK}{\mathscr{K}}
\newcommand{\eM}{\mathscr{M}}
\newcommand{\eU}{\mathscr{U}}
\newcommand{\eV}{\mathscr{V}}
\newcommand{\fC}{\mathfrak{C}}
\newcommand{\fG}{\mathfrak{G}}
\newcommand{\fI}{\mathfrak{I}}
\newcommand{\fM}{\mathfrak{M}}
\newcommand{\fU}{\mathfrak{U}}
\newcommand{\fX}{\mathfrak{X}}

\newcommand{\bA}{\mathbb{A}}
\newcommand{\bC}{\mathbb{C}}
\newcommand{\bD}{\mathbb{D}}
\newcommand{\bP}{\mathbb{P}}
\newcommand{\bQ}{\mathbb{Q}}
\newcommand{\bZ}{\mathbb{Z}}
\newcommand{\sC}{\mathcal{C}}
\newcommand{\sE}{\mathcal{E}}
\newcommand{\sL}{\mathcal{L}}
\newcommand{\sO}{\mathcal{O}}
\newcommand{\sT}{\mathcal{T}}
\newcommand{\age}{\mathrm{age}}
\newcommand{\ch}[2]{\genfrac{(}{)}{0pt}{}{#1}{#2}}
\newcommand{\edim}{\mathrm{edim}}
\newcommand{\pd}{\mathbb{P}^2_{D,2}}
\newcommand{\rk}{\mathrm{rk}}
\newcommand{\tw}{\mathrm{tw}}

\begin{document}

\def\baselinestretch{1}
\def\arraystretch{1.6}

\begin{abstract}
The Gromov-Witten theory of Deligne-Mumford stacks is a recent development, and hardly any computations have been done beyond $3$-point genus $0$ invariants.  This paper provides explicit recursions which, together with some invariants computed by hand, determine all genus $0$ invariants of the stack $\bP^2_{D,2}$.  Here $D$ is a smooth plane curve and $\bP^2_{D,2}$ is locally isomorphic to the stack quotient $[U/(\bZ/(2))]$, where $U\to V\subset \bP^2$ is a double cover branched along $D\cap V$.  The introduction discusses an enumerative application of these invariants.
\end{abstract}

\maketitle

\section{Introduction}
\label{sec:intro}
\setcounter{equation}{0}

This paper provides a family of examples where the associativity relations together with a small number of invariants determine all genus $0$ Gromov-Witten invariants of a Deligne-Mumford stack.  It is by now a standard technique to use these relations to compute genus $0$ invariants of smooth projective varieties.  This is partly due to the first reconstruction theorem, which says that if the algebraic cohomology ring of a smooth variety is generated by divisors, then its genus $0$ Gromov-Witten invariants are determined by the $n$-point invariants for $n\le 2$ \cite{KM}.  It is an interesting question whether any kind of reconstruction theorem holds for stacks.

The stacks we study are parametrized by smooth plane curves.  If the degree of the curve $D$ happens to be even, then $\pd$ is the stack quotient of a double cover of $\bP^2$ which is branched along $D$ by the Galois group of the cover.  If the degree of $D$ is odd, then this picture only holds locally, but the quotient stacks glue together.  The general construction of $\pd$ is worked out in \cite{stacks}.

The Gromov-Witten invariants take as input numerical equivalence classes on the inertia stack of $\pd$.  In this paper we instead use the coarse moduli space of the inertia stack, which is isomorphic to a disjoint union of $\bP^2$ with $D$.  The interesting classes are the class of a point in $\bP^2$, the fundamental class of $D$, and the class of a point in $D$, and we denote them $T_2$, $T_3$, and $T_4$ respectively.

A twisted stable map is a representable morphism of stacks $\fC\to\pd$ whose induced morphism of coarse moduli spaces $C\to\bP^2$ is an ordinary stable map such that all ``twisted'' points of $\fC$ are either nodes or marked points.  A twisted point is a point which has a nontrivial stabilizer group, and for $\pd$ this group will always be $\bZ/(2)$.  The stack $\fC$ is required to have a certain local picture at nodes and marked points \cite[\S 2.1]{ACV}, and this local picture completely determines $\fC$ given the set of twisted points and the orders of their cyclic stabilizer groups.

Twisted stable maps to $\pd$ are best understood in terms of tangencies.  A morphism $f:C\to\bP^2$ from a smooth curve $C$ such that $f^{-1}(D)$ is zero dimensional lifts to a morphism $f:C\to\pd$ if and only if $f^*(D)=2E$ for some divisor $E$ \cite[3.3.6]{stacks}.  So we can loosely say that a twisted stable map $F:\fC\to\pd$ has to have an even order contact with $D$ at every untwisted point of $\fC$.  If $f:C\to\bP^2$ is the coarse map of $F$, then the order of contact between $C$ and $D$ at a twisted point of $\fC$ must be odd [ibid].  At a twisted or untwisted node, this condition must hold for both components on which the nodal point lies \cite{enumerate}.

This work was motivated by the idea that Gromov-Witten invariants of $P^2_{D,2}$ might count rational plane curves having prescribed contacts with $D$.  While this is often true, there are plenty of counterexamples.  When $D$ is a line, the numbers of such curves can be computed from a recursion of Caporaso and Harris \cite[\S 1.4]{CH}.  The first discrepancy between the Gromov-Witten invariants and the actual numbers occurs in degree $4$ with the invariant $I_4(T_2^7T_4^4)$.  If this invariant were enumerative, it would count the number of rational quartics passing through $7$ general points in $\bP^2$ and $4$ general points on $D$.  However, the actual number is $398$, while the invariant is $416$.  This is explained by the existence of twisted stable maps from reducible curves into $\pd$ which are illustrated by the diagram below.  The dots correspond to twisted points and the numbers indicate the degree of the component they label.
\medskip

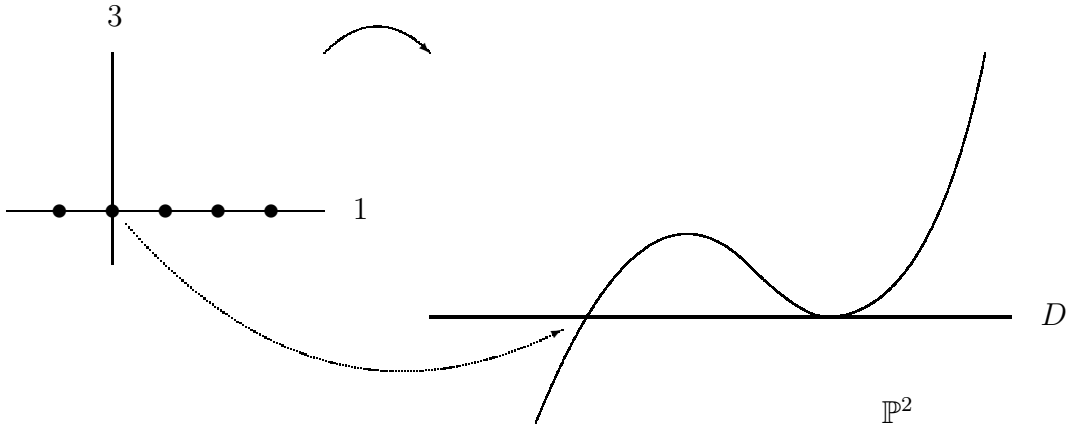
\begin{picture}(450,170)(0,30)

\put(60,100){\line(0,1){80}}
\put(57,190){$3$}
\put(20,120){\line(1,0){120}}
\put(150,117){$1$}
\multiput(40,120)(20,0){5}{\circle*{5}}

\thicklines
\put(180,80){\line(1,0){220}}
\thinlines
\put(410,77){$D$}
\qbezier(390,180)(370,80)(330,80)
\qbezier(330,80)(320,80)(300,100)
\qbezier(300,100)(260,140)(220,40)

\put(350,40){$\bP^2$}

\qbezier[75](140,180)(160,200)(180,180)
\put(180,180){\vector(1,-1){0}}

\qbezier[150](65,115)(140,30)(230,75)
\put(230,75){\vector(2,1){0}}

\end{picture}

This shows that there is a contribution to the Gromov-Witten invariant coming from rational cubics which are tangent to the line and pass through the $7$ general points.  There are $36$ such curves, but each of them counts with multiplicity $1/2$, because the twisted curve in the diagram has a ``ghost automorphism'' which 2-commutes with the map.  Ghost automorphisms are automorphisms of the twisted curve which cover the identity morphism on its coarse moduli space, and it is shown in \cite[7.1.1]{ACV} that they all come from twisted nodes.  It is important that the component mapping onto $D$ has an odd number of twisted points, because otherwise a representable morphism would not exist.

We have compared the Gromov-Witten invariants $$I_d(T_2^{3d-1-a-b}T_3^{d-a-2b}T_4^{a})$$ with the number of rational degree $d$ curves $C$ which have $b$ tangencies with $D$ at smooth points of $C$, meet $D$ transversely at $a$ general points, and meet $3d-1-a-b$ general points in $\bP^2$.  For all except one case with $d\le 6$, the invariants either give the same answer as the recursion of Caporaso and Harris or else the discrepancy can be accounted for in a manner similar to the above example.  The remaining case, $I_6(T_2^{11}T_4^6)$, involves a virtual fundamental class computation which we have not done.

When $D$ is a conic, there are similar contributions which prevent the invariants from being enumerative, but when $\deg(D)\ge 3$ no rational curve can map onto $D$, so such contributions are eliminated.  This does not mean that they become enumerative, however.  When $D$ is a quartic, there are components consisting entirely of multiple covers which contribute to the invariants.  The simplest example is $I_2(T_2)$, which if enumerative would count conics passing through a point which are tangent to $D$ four times.  The extra contribution comes from double covers of lines which are tangent to $D$ once such that the branch points of the cover lie at the transverse intersection points.  This type of contribution can occur whenever $\deg(D)\ge 4$.

When $D$ is a cubic, the moduli space of twisted stable maps still has components consisting of multiple covers, but they do not contribute to the Gromov-Witten invariants.  In \cite{enumerate}, it is shown that all of the positive degree genus $0$ Gromov-Witten invariants of $\pd$ are enumerative.  Stronger enumerative results have recently been obtained.  In \cite{enum2}, more general ``$r$th root'' stacks were used to compute the numbers of rational degree $d$ curves with arbitrary contact conditions imposed relative to a smooth cubic, except for the case of a single contact of order $3d$.

By this point, the reader may have guessed that the genus $0$ Gromov-Witten invariants of $\pd$ do not at all behave like those of homogeneous spaces.  The moduli stack of stable maps to $\pd$ is generally singular and fails to have the expected dimension.  Even in degree $0$ the virtual fundamental class plays a role, and a surprising consequence is the existence of $n$-pointed degree $0$ invariants with $n>3$.  For stacks in general, such invariants can exist for $n$ arbitrarily large, but for $\pd$ we are fortunate that the only such invariant occurs with $n=4$.  Moreover, this extra invariant can be computed from the associativity relations.

When $D$ is a line or a conic, there exist infinitely many invariants of degree $\delta:=\deg(D)$.  These come from stable maps which on the level of coarse moduli spaces have image $D$, which is why they only occur when $D$ is rational.  If $$\lambda_k = \frac{(-1)^k k !}{2^{k+1}} \mbox{ for } k\ge 0,$$ then when $D$ is a line, the Gromov-Witten invariant $I_1(T_3^kT_4^{k+3})$ equals $\lambda_k$, while if $D$ is a conic, $I_2(T_3^kT_4^{k+6})$ equals $\lambda_k$.

The paper is organized as follows.  Sections \ref{sec:maps} and \ref{sec:vfc} contain general results about stable maps into $X_{D,r}$.  Section \ref{sec:maps} defines evaluation maps and Gromov-Witten invariants, while section \ref{sec:vfc} deals with the virtual fundamental class and works out the expected dimension of the stack of stable maps.  Along the way, we investigate the tangent bundle of $X_{D,r}$ in section \ref{sec:vfc}.3, showing that it is related to the tangent bundle of $X$ by elementary transformations.  Section \ref{sec:computations} contains direct computations of some degree 0 and 1 invariants which are needed in order to use the recursions.  Finally, in section \ref{sec:quant_prod} we define the big quantum product, work out the recursions, and provide an algorithm which can be used to compute any genus 0 invariant.
\smallskip

\noindent{\bf Notation and conventions}

Throughout this paper, all schemes are equipped with structure morphisms to $\spec\bC$ and all morphisms of schemes respect these structure morphisms.  As a consequence, the same is true for all stacks that appear.  We use $\mu_r$ to denote the group of $r$th roots of unity in $\bC$.  Given a stack $\fX$, we use $\widehat{\fX}$ to denote the coarse moduli space of its inertia stack, assuming it exists.

In this paper, we deviate slightly from \cite{AGV} where Gromov-Witten invariants are defined in terms of a stack of twisted stable maps with sections of all gerbes, which is denoted $\bar{\eM}_{g,n}(\fX,\beta)$.  We prefer instead to work with the stack $\eK_{g,n}(\fX,\beta)$, which is the one defined in \cite{AV}.  We implicitly assume that all twisted stable maps are balanced.  In addition, we work with evaluation morphisms from $\eK_{g,n}(\fX,\beta)$ to the coarse moduli space of the inertia stack, rather than the rigidification of the inertia stack used in \cite{AGV}.
\smallskip

\noindent{\bf Acknowledgements}

This paper was derived from part of my Ph.D. thesis at Columbia University under the supervision of Michael Thaddeus.  I would like to thank Dan Abramovich and Martin Olsson for helpful discussions.

\section{Stable maps to \boldmath $r$th root stacks}
\label{sec:maps}

Let $X$ be a smooth projective variety over $\bC$, $D\subseteq X$ be a smooth divisor, and $r$ be a positive integer.  In this section, we investigate stable maps into $X_{D,r}$ and define Gromov-Witten invariants, though we leave a discussion of the virtual fundamental class for later.  First we introduce the contact type of a twisted stable map and define evaluation maps in terms of this.  It is straightforward to show that the definition in \cite{AGV} gives rise to the same evaluation maps.

\subsection{The $r$th root construction}

Here is a summary of results from section 2 of \cite{stacks} which are needed in this paper.  Given a scheme $X$, an effective Cartier divisor $D\subseteq X$, and a positive integer $r$ which is invertible on $X$, there is a Deligne-Mumford stack $X_{D,r}$ over $X$ on which the pair $(\sO_X(D),s_D)$ has an $r$th root.  More precisely, if $f:S\to X$ is a morphism of schemes, then an object of $X_{D,r}$ over $f$ consists of a triple $(M,t,\varphi)$, where $M$ is an invertible sheaf on $S$, $t$ is a global section of $M$, and $\varphi:M^r\to f^*\sO_X(D)$ is an isomorphism such that $\varphi(t^r)=f^*s_D$.  There is a universal object on $X_{D,r}$ covering the projection $\pi:X_{D,r}\to X$.  We denote the universal line bundle $\sT$ and its section $\tau$.  The vanishing locus of $\tau$ is a substack $\fG\subseteq X_{D,r}$ which is mapped into $D$ by $\pi$, and the restriction $\fG\to D$ is an \'etale gerbe with band $\mu_r$.  The projection $\pi$ is an isomorphism away from $D$, and it is ramifed over $D$.  In fact, $\pi^{-1}(D)$ is the $r$th order infinitesimal neighborhood of $\fG$ in $X_{D,r}$.

\subsection{Twisted curves}

Recall that a prestable curve over a Noetherian scheme $S$ is a flat, proper morphism $C\to S$ whose geometric fibers are curves with at worst nodal singularities.  An $n$-marked prestable curve has in addition $n$ sections $s_i:S\to C$ whose images are disjoint and don't intersect the singular locus of $C\to S$.  To give $n$ such sections is equivalent to giving $n$ disjoint effective Cartier divisors $D_i\subseteq C$ which map isomorphically to $S$ \cite[\S 4]{stacks}.  The passage from curves to twisted curves requires the addition of stack structure at the nodes and markings \cite[4.1.2]{AV}.  The stack structure at a marking is always obtained from the coarse moduli space by applying the $r$th root construction to the image of the section.  We do not need to consider the stack structure at the node, but for an interesting treatment see \cite{Ol}.

Given a Deligne-Mumford stack $\fC$, $n$ divisors $D_i$ and $n$ positive integers $r_i$, we use $\fC_{\bD,\vec{r}}$ for the result of applying the $r_1$th root construction along $D_1$, followed by the $r_2$th root construction along $D_2$, and so on.

\begin{prop}
\label{twisted_curve}
To give a twisted nodal $n$-pointed genus $g$ curve over a connected Noetherian scheme $S$ is equivalent to giving
\begin{enumerate}
  \item a nodal $n$-pointed genus $g$ curve $C$ over $S$ with markings $D_1,\ldots, D_n$,
  \item for each marking $D_i\subseteq S$, a positive integer $r_i$,
  \item and a twisted nodal $0$-pointed genus $g$ curve $\fC$ over $S$, whose coarse moduli scheme is $C$.
\end{enumerate}
Then the twisted curve is isomorphic to $\fC_{\bD,\vec{r}}$ in such a way that the $i$th marking of the twisted curve is sent to the gerbe over $D_i$.
\end{prop}

This can be proven in the same way that \cite[4.1]{stacks} was, and also follows from \cite{Ol}.  Note that $\fC_{\bD,\vec{r}}\cong \fC\times_C C_{\bD,\vec{r}}$.  Moreover, the divisors $D_i$ don't pass through any nodes of the fibers, so the $r_i$th root construction along $D_i$ is taking place away from the points of $\fC$ which have stack structure.  While this statement is not very precise, it is meant to assure the reader that much of what was done in \cite{stacks} carries over without change.  In particular, the notion of contact type extends to arbitrary stable maps $\fC_{\bD,\vec{r}}\to X_{D,r}$.

\subsection{Contact type}

In \cite[\S 3.2]{stacks}, we derived some key results about line bundles and sections on $C_{\bD,\vec{r}}$.  They extend without change to $\fC_{\bD,\vec{r}}$, because the results of [ibid., \S 3.1], from which they follow, are stated in sufficient generality.  We state the more general results here.

Let $S$ be a \emph{connected}, Noetherian scheme and let $C$, $\fC$, and $D_i$ be as in Proposition \ref{twisted_curve}.  Let $\gamma:\fC_{\bD,\vec{r}}\to\fC$ be the projection, and let $\sT_i$ be the tautological sheaf associated to the $r_i$th root construction along $D_i$.  The important properties of these sheaves are contained in the next two corollaries.

\begin{cor}
\label{inv_sheaf}
Let $\sL$ be an invertible sheaf on $\fC_{\bD,\vec{r}}$.  Then there exist an invertible sheaf $L$ on $\fC$ and integers $k_i$ satisfying $0\le k_i\le r_i-1$ such that
$$\sL\cong\gamma^*L\otimes\prod_{i=1}^n\sT_i^{k_i}.$$
Moreover, the integers $k_i$ are unique, $L$ is unique up to isomorphism, and $\sT_i^{r_i}\cong\gamma^*\sO_{\fC}(D_i)$.
\end{cor}

\begin{cor}
\label{inv_sheaf_section}
Given the decomposition in the previous corollary, every global section of $\sL$ is of the form $\gamma^*s\otimes\tau_1^{k_1}\otimes\cdots\otimes\tau_n^{k_n}$ for a unique global section $s$ of $L$, where $\tau_i$ is the tautological section of $\sT_i$.
\end{cor}

If $F:\fC_{\bD,\vec{r}}\to X_{D,r}$ is any morphism, then Corollary \ref{inv_sheaf} associates to $F^*\sT$ a unique $n$-tuple of integers $k_1,\ldots,k_n$, where $\sT$ is the tautological sheaf of $X_{D,r}$.  Let $U\subseteq C$ be the complement of the divisors $D_i$ and let $\fU\subseteq\fC$ be its preimage.  Note that $\fU$ is isomorphic to its preimage in $\fC_{\bD,\vec{r}}$.

\begin{prop}
\label{rep_criterion}
The morphism $F$ is representable if and only if its restriction to $\fU$ is representable and for every $i$, $r_i$ divides $r$ and $k_i$ is relatively prime to $r_i$.
\end{prop}

The proof goes through as in \cite[3.3.3]{stacks}.  Being representable is a pointwise condition by \cite[4.4.3]{AV}, so for $F$ to be representable is equivalent to its restrictions to $\fU$ and the gerbes over $D_1,\ldots,D_n$ to be representable.

Now assume that $F$ is representable.

\begin{definition}
\label{contact_type}
Let $\varrho_i=k_ir/r_i$ for $1\le i\le n$.  We define the \emph{contact type} of $F$ to be the $n$-tuple $\vec{\varrho}=(\varrho_1,\ldots,\varrho_n)$.
\end{definition}

Multiplying the identity $\gcd(k_i,r_i)=1$ by $r/r_i$ yields $\gcd(\varrho_i,r)=r/r_i=\varrho_i/k_i$.  Thus $\vec{k}$ and $\vec{r}$ are determined by $\vec{\varrho}$ via the following formulas.
\begin{equation}
\label{rk_dependence}
r_i=\frac{r}{\gcd(\varrho_i,r)}\;\;\;\;\;\;
k_i=\frac{\varrho_i}{\gcd(\varrho_i,r)}
\end{equation}

For any $\beta\in N_1(X)$, the contact type determines $n$ locally constant functions on $\eK_{g,n}(X_{D,r},\beta)$ with integer values from $0$ to $r-1$.  We define $\eK_{g,n}(X_{D,r},\beta,\vec{\varrho})$ to be the open and closed substack of $\eK_{g,n}(X_{D,r},\beta)$ consisting of stable maps with contact type $\vec{\varrho}$.  In \cite{stacks}, we found that a stable map from a smooth twisted curve into $X_{D,r}$ which does not map into the gerbe has to have a contact type for which $D\cdot\beta-\sum\varrho_i$ is a multiple of $r$.  Equation \ref{eq:edim} implies that this holds in general.

\subsection{Evaluation maps}

In \cite{AGV}, evaluation maps are defined which go from $\eK_{g,n}(X_{D,r},\beta)$ to a rigidification of the inertia stack of $X_{D,r}$.  This allows for a richer structure than is required for our purposes.  We instead define evaluation maps which go to the coarse moduli scheme of the inertia stack, which we denote $\widehat{X_{D,r}}$.  It turns out that $\widehat{X_{D,r}}$ is isomorphic to a disjoint union of $X$ with $r-1$ copies of $D$.  We denote these copies of $D$ by $(\widehat{X_{D,r}})_i$ for $1\le i\le r-1$, and let $(\widehat{X_{D,r}})_0$ denote the component $X$.  A morphism from a connected scheme $S$ to $\widehat{X_{D,r}}$ is therefore determined by a morphism $S\to X$ and an integer from $0$ to $r-1$ which is positive only if $S\to X$ factors through $D$.

Given an object $\xi$ of $\eK_{g,n}(X_{D,r},\beta)$ over a connected Noetherian scheme $S$, let $e_i(\xi)$ be the pair $(g,\varrho_i)$, where $\varrho_i$ is the $i$th component of the contact type of $\xi$ and $g:S\to X$ is the composition of the $i$th section $S\to D_i\subseteq\fC$ with the morphism $\fC\to X$.
$$\xymatrix{
 & \fC_{\bD,\vec{r}}\ar[r]^F\ar[d]_{\gamma} & X_{D,r}\ar[d] \\
D_i \ar@{^{(}->}[r] & \fC\ar[r]^f & X}$$

In light of the following lemma, this defines a morphism $e_i:\eK_{g,n}(X_{D,r},\beta)\to \widehat{X_{D,r}}$.

\begin{lemma}
\label{ev_map}
If $\varrho_i\ne 0$, then the morphism $g:S\to X$ factors through $D$.
\end{lemma}

\begin{proof} The morphism $F$ determines an invertible sheaf $M$ on $\fC_{\bD,\vec{r}}$, which by Corollary \ref{inv_sheaf} is isomorphic to $\gamma^*L\otimes\prod_{j=1}^n\sT_j^{k_j}$.  It also determines a section $t$ of $M$, which by Corollary \ref{inv_sheaf_section} corresponds to $\gamma^*s\prod \tau_j^{k_j}$.  Moreover $M^r$ is isomorphic to $\gamma^*f^*\sO_X(D)$ in such a way that $t^r$ goes to $\gamma^*f^*s_D$.  Taking $r$th powers of our expressions for $M$ and $t$, the uniqueness assertions of Corollaries \ref{inv_sheaf} and \ref{inv_sheaf_section} imply that $L^r(\sum\varrho_j D_j)\cong f^*\sO_X(D)$ by an isomorphism sending $s^r\prod s_{D_j}^{\varrho_j}$ to $f^*s_D$.  It follows that if $\varrho_i\ne 0$, then $f^*s_D$ vanishes on $D_i$, which proves the lemma.\end{proof}

\subsection{Gromov-Witten invariants}

To define Gromov-Witten invariants, one must first define the virtual fundamental class $[\eK_{g,n}(X_{D,r},\beta)]^v$.  This is an element of the Chow group of $\eK_{g,n}(X_{D,r},\beta)$ which has the expected dimension.  For now we assume it has already been defined.

Fix integers $g,n\ge 0$ and an effective class $\beta\in N_1(X)$.  Let $\vec{e}=e_1\times\cdots\times e_n:\eK_{g,n}(X_{D,r},\beta)\to \widehat{X_{D,r}}^{\times n}$ and let $p:\widehat{X_{D,r}}^{\times n}\to\spec\bC$ be the structure morphism.  Since $\eK_{g,n}(X_{D,r},\beta)$ is proper over $\bC$ \cite[1.4.1]{AV}, there is a proper pushforward $\vec{e}_*$ on Chow groups with rational coefficients defined by \cite[3.6]{Vi}.

Given any numerical equivalence classes $a_1,\ldots,a_n\in N^*_{\bQ}(\widehat{X_{D,r}})$, the Gromov-Witten invariant is defined to be
\begin{equation}
\label{def:GWI}
I^{g}_{\beta}(a_1\cdots a_n):= p_*\left((a_1\times\cdots\times a_n) \cap \vec{e}_*[\eK_{g,n}(X_{D,r},\beta)]^v\right),
\end{equation}
which is regarded as an element of $\bQ$.  As the definition is clearly symmetric in $a_1,\ldots,a_n$, we regard $I^g_{\beta}$ as an element of the dual of $\Sym^*(N^*_{\bQ}(\widehat{X_{D,r}}))$.  The following proposition is an immediate consequence of the definitions.

\begin{prop}
\label{invts_zero}
Suppose $a_i=[V_i]$ for some (irreducible) subvariety $V_i\subseteq \widehat{X_{D,r}}$ for $1\le i\le n$.  Let $\vec{\varrho}$ be such that $V_i\subseteq (\widehat{X_{D,r}})_{\varrho_i}$ for each $i$.  If there are no stable maps $\fC\to X_{D,r}$ in $\eK_{g,n}(X_{D,r},\beta,\vec{\varrho})$ such that the $i$th marked point is sent under $\fC\to X_{D,r}\to X$ into the image of $V_i\subseteq \widehat{X_{D,r}}\to X$ for all $i$, then $I_{\beta}^g(a_1\cdots a_n)=0$.
\end{prop}

\section{Virtual Fundamental Class}
\label{sec:vfc}

\subsection{Preliminaries}

This section is devoted to studying the virtual fundamental class for the moduli space of stable maps into $X_{D,r}$.  To begin, we consider an arbitrary smooth, proper, Deligne-Mumford stack $\fX$ which has a projective coarse moduli scheme.  Note that $X_{D,r}$ is smooth by Proposition \ref{smoothness}.  The reader may also wish to consult \cite[\S 4.6]{AGV}.

The virtual fundamental class of $\eK:=\eK_{g,n}(\fX,\beta)$ is constructed in the same way as it is for varieties \cite{Be,BF,LT}.  Here we adopt the approach of Behrend and Fantechi.  Let $\fM:=\fM_{g,n}^{\tw}$ be the smooth Artin stack of twisted curves with no stability condition.  This is obtained from \cite[1.9]{Ol} by base extension to $\spec\bC$.  We will define an element $E^{\bullet}$ of the derived category of $\eK$ which is equivalent to a two-term complex of locally free sheaves in positions $-1$ and $0$.  Then we will define a morphism $\phi:E^{\bullet}\to L_{\eK/\fM}$ which is surjective on $h^{-1}$ and an isomorphism on $h^0$.  After verifying that $\phi$ is a perfect relative obstruction theory, the construction of \cite[\S 7]{BF} produces a virtual fundamental class $[\eK]^v\in A_*(\eK)$.  Note that the mapping cone construction allows one to go from a perfect relative obstruction theory to an absolute one which determines the same virtual fundamental class \cite[App. B]{GP}.

Before proceeding to define $\phi$, we state some consequences of this construction.  On any connected component of $\eK$, $[\eK]^v$ is homogeneous of degree $\rk(E^{\bullet})+\dim(\fM)$, which is called the expected dimension.  By \cite[\S 4]{BF}, a perfect obstruction theory for $\eK$ gives rise to an obstruction theory in the classical sense at any point of $\eK$.   It follows that at every point of $\eK$, the dimension is greater than or equal to the expected dimension with equality implying that $\eK$ is a local complete intersection at that point \cite[7.4]{Ha2}.

Suppose that $\eU\subseteq\eK_{g,n}(\fX,\beta)$ is an open substack which has the expected dimension.  By \cite[5.10]{BF}, the restriction of the virtual fundamental class to $\eU$ is equal to the virtual fundamental class determined by the restriction of $\phi$ to $\eU$.  Since $\eU$ has the expected dimension, the kernel of $$h^{-1}(\phi\vert_{\eU}):h^{-1}(E^{\bullet}\vert_{\eU}) \to h^{-1}(L_{\eU/\fM})$$ has rank $0$.  Moreover, it is torsion free since it is a subsheaf of a locally free sheaf.  Therefore, it equals $0$ and $\phi\vert_{\eU}$ is an isomorphism.  The following proposition follows by definition of the virtual fundamental class \cite[\S 5]{BF}.

\begin{prop}
\label{vfc_coeff=1}
If $\eU\subseteq\eK_{g,n}(\fX,\beta)$ is an open substack which has the expected dimension, then the restriction of $[\eK_{g,n}(\fX,\beta)]^v$ to $\eU$ is $[\eU]$, the ordinary fundamental class.
\end{prop}

\subsection{Perfect relative obstruction theory}

We continue with our previous notation.  Let $\eC$ be the universal curve over $\eK$.
$$\xymatrix{
\eC \ar[r]^{f} \ar[d]_{\pi} & \fX \\
\eK}$$
There exist canonical log structures on $\eC$ and $\eK$ which make $\pi$ log smooth \cite[\S 3]{Ol}.  Let $\omega$ be the relative sheaf of log differentials, so that $\omega$ is a dualizing sheaf for $\pi$.

Let $$E^{\bullet}=(R\pi_*(f^*T_{\fX}))^{\vee}.$$  By duality, $E^{\bullet}$ is isomorphic to $R\pi_*(f^*L_{\fX}\otimes\omega)$.  Let $\eC^{\tw}\to\fM$ be the universal twisted curve.  We have morphisms $f^*L_{\fX}\to L_{\eC} \to L_{\eC/\eC^{\tw}}$ and an isomorphism $\pi^*L_{\eK/\fM}\to L_{\eC/\eC^{\tw}}$, the latter coming from the following fiber square.
$$\xymatrix{
\eC \ar[r] \ar[d]_{\pi} \ar@{}[dr]|{\Box} & \eC^{\tw} \ar[d] \\
\eK \ar[r] & \fM}$$
Composing, we obtain a morphism $f^*L_{\fX}\to\pi^*L_{\eK/\fM}$, which by duality is equivalent to a morphism $\phi:E^{\bullet}\to L_{\eK/\fM}$.

We show that $E^{\bullet}$ is equivalent to a two-term complex of locally free sheaves (it is enough for this to hold locally on $\eK$, but it is true globally).  Let $\overline{\pi}:\overline{\eC}\to\eK$ be the relative coarse moduli space of $\eC$ over $\eK$, i.e. the fibers of $\overline{\pi}$ are the coarse moduli spaces of the fibers of $\pi$.  Let $\gamma:\eC\to\overline{\eC}$ be the projection.  By \cite[4.2]{OS}, $\gamma_*f^*T_{\fX}$ is a coherent sheaf which is flat over $\eK$, and $R^i\gamma_*f^*T_{\fX}=0$ for $i>0$.  It follows that $E^{\bullet}\cong R\overline{\pi}_*(\gamma_*f^*T_{\fX})$.  One can now apply \cite[5]{Be}, since the proof of that proposition requires only that $E$ is a coherent sheaf, flat over $T$.

Since deformation theory for Deligne-Mumford stacks works in the same way as for schemes (all the results of \cite{Il} apply), the remainder of the proof that $\phi$ is a perfect relative obstruction theory is identical to the case where $\fX$ is a smooth projective variety.

In summary, if $\sC=(f:\fC\to\fX,\Sigma_1,\ldots,\Sigma_n)$ is a stable map in $\eK_{g,n}(\fX,\beta)$, then the expected dimension of $\eK_{g,n}(\fX,\beta)$ at $\sC$ is
\begin{equation}
\label{eq:edim_gen}
\chi(f^*T_{\fX})+\dim(\fM_{g,n}^{\tw}).
\end{equation}

\subsection{Riemann-Roch for twisted curves}

Here we recall a version of the Riemann-Roch theorem for balanced nodal twisted curves which is proven in \cite{AGV2}.  It will be applied later to $f^*TX_{D,r}$, where $f:\fC\to X_{D,r}$ is a twisted stable map.  First we recall the definition of age.

Let $\fC$ be a balanced twisted nodal curve, and let $\sE$ be a locally free sheaf on $\fC$.  For any closed point $x\in\fC$ there is a closed substack $\fG_x\subseteq\fC$ called the residue gerbe of $x$ \cite[\S 11]{LMB}.  It is isomorphic to $B\mu_r$ for some positive integer $r$ called the \emph{order of twisting} of $x$.  If $r>1$, $x$ is called a \emph{twisted point}, and otherwise it is called \emph{untwisted}.  The restriction of $\sE$ to $B\mu_r$ determines a representation of $\mu_r$, but this is not canonical since it depends on the isomorphism $B\mu_r\to\fG_x$.  If $x$ is a smooth point, then we can choose an isomorphism $B\mu_r\to\fG_x$ so that the restriction of the tangent bundle $T\fC$ to $B\mu_r$ corresponds to the standard representation of $\mu_r$.  Then the representation of $\mu_r$ determined by $\sE$ is well-defined, and we denote this representation by $\sE_x$.

At nodes, one can study $\sE$ by pulling back to the normalization of $\fC$ \cite[1.18]{Vi}.  Locally in the \'etale topology, the node looks like the stack quotient of $\spec\bC[x,y]/(xy)$ by the group $\mu_r$ which acts by $\zeta\cdot (x,y) = (\zeta x,\zeta^{-1}y)$.  The normalization of this stack is two copies of the affine line each with the origin twisted to order $r$.  So a nodal point $x\in\fC$ has two preimages $x_i$ in the normalization of $\fC$, and each determines canonically a representation $\sE_{x_i}$ of $\mu_r$.  The fact that $\fC$ is balanced (encoded in the action above) implies that these are dual representations.

\begin{definition}
\label{age_def}
The \emph{age} of a locally free sheaf $\sE$ at a smooth point $x$ of a twisted curve $\fC$ is a rational number denoted $\age(\sE,x)$.  Let $V$ denote the standard representation of $\mu_r$, where $r$ is the order of twisting of $x$.  If $$\sE_x\cong\oplus_{i=1}^{n} V^{\otimes k_i}$$ where $0\le k_i\le r-1$, then $$\age(\sE,x)=\sum_{i=1}^n \frac{k_i}{r}.$$
\end{definition}

\begin{thm}
\label{Riemann-Roch}
Let $\fC$ be a balanced twisted nodal curve of genus $g$ over $\bC$ and let $\sE$ be a locally free sheaf on $\fC$.  Then $$h^0(\fC,\sE)-h^1(\fC,\sE)=\deg(\sE)+\rk(\sE)(1-g)-\sum_x \age(\sE,x),$$ where the sum is over the twisted smooth points of $\fC$.
\end{thm}

\subsection{Tangent bundle of \boldmath $X_{D,r}$}
\label{subsec:tbundle}

To apply the Riemann-Roch theorem to $f^*TX_{D,r}$, we need some results about the tangent bundle to $X_{D,r}$.  Let $\sT$ be the tautological sheaf on $X_{D,r}$ (see 2.1) and let $\fG$ be the gerbe of $X_{D,r}$.  We use $\fG_{k}$ to denote the $k$th order infinitesimal neighborhood of $\fG$.  Let $\pi:X_{D,r}\to X$ be the projection.

Recall from \cite[2.4.1]{stacks} that $X_{D,r}$ is locally a stack quotient in the following way.  If $V=\spec S\subseteq X$ is an affine open set on which we have a trivialization of $\sO_X(D)$, then the tautological section of $\sO_X(D)$ corresponds to an element $\alpha\in S$, and if $U=\spec A$, where $A=S[x]/(x^r-\alpha)$, then $X_{D,r}\times_X V$ is isomorphic to $[U/\mu_r]$, where $\mu_r$ acts on $A$ by $(\zeta,x)\mapsto \zeta^{-1}x$ and $(\zeta,s)\mapsto s$ for $\zeta\in\mu_r$ and $s\in S$.  Moreover, the preimage of $\fG_k$ in $U$ is defined by the vanishing of $x^k$.  This shows that $\fG_r=\pi^{-1}(D)$.

We assume that $X$ is a smooth projective variety and $D\subseteq X$ is a smooth divisor.  First we show that $X_{D,r}$ is smooth, so that $TX_{D,r}$ is locally free.

\begin{prop}
\label{smoothness}
$X_{D,r}$ is smooth over $\bC$.
\end{prop}

\begin{proof} Using the above notation, it suffices to show that $U$ is smooth.  Since $U$ is a closed subscheme of the smooth variety $V\times\bA^1$ defined by $x^r-\alpha=0$, it suffices to show that for any closed point $p\in V\times\bA^1$, $x^r-\alpha\notin m_p^2$, where $m_p\subseteq\sO_{V\times\bA^1,p}$ is the maximal ideal.

If $x^r-\alpha\in m_p^2$, then $\partial/\partial x (x^r-\alpha)\in m_p$, so $rx^{r-1}\in m_p$, which implies $x\in m_p$.  It follows that $p\in D\times\{0\}$ and $\alpha\in m_p^2$, which contradicts the smoothness of $D$.\end{proof}

\begin{prop}
\label{tan_exact}
There is an exact sequence of sheaves on $X_{D,r}$:
$$0\to TX_{D,r}\to\pi^*TX\to\pi^*\sO_X(D)\otimes\sO_{\fG_{r-1}}\to 0.$$
\end{prop}

\begin{proof} Since $\pi^*\sO_D\cong\sO_{\fG_r}$, the morphisms in the above sequence clearly exist.  The injectivity of the differential follows from the fact that it is an isomorphism away from $\fG$.  Note also that
$$\pi^*TX\to\pi^*(TX\otimes\sO_D)\to\pi^*(\sO(D)\otimes\sO_D)\to\pi^*\sO(D)\otimes\sO_{\fG_{r-1}}$$
is a sequence of surjective maps.  It remains to show exactness in the middle, and this can be shown locally.  It suffices to show this after pulling back to $U$ since $U\to X_{D,r}$ is \'etale.  So we need to show that the following sequence of $A$-modules is exact in the middle.
$$0\to\Der_{\bC}(A,A)\to\Der_{\bC}(S,S)\otimes_S A\to \alpha^{-1}A/x^{-1}A \to 0$$
The second homomorphism is defined by sending a derivation $\delta:S\to S$ to $\alpha^{-1}\delta(\alpha)$.  The composition is $0$ because $\alpha^{-1}\delta(x^r)=rx^{-1}\delta(x)$.  Note that a derivation $A\to A$ is uniquely determined by a derivation $\delta:S\to A$ and an element $\delta(x)\in A$ such that $rx^{r-1}\delta(x)=\delta(\alpha)$.  If $\sum a_i\delta_i\mapsto 0$ in $\alpha^{-1}A/x^{-1}A$, then $\sum a_i\delta_i(\alpha)\in\alpha x^{-1}A=x^{r-1}A$, so $\sum a_i\delta_i$ extends to a derivation $A\to A$.\end{proof}

Let $\sE$ be the kernel of $TX\to\sO_D(D)$, which is locally free.  From the above proposition, it follows that as subsheaves of $\pi^*TX$, $\pi^*\sE$ is a subsheaf of $TX_{D,r}$.  We have the following commutative diagrams with exact rows.
$$\xymatrix{
0\ar[r] & \pi^*\sE \ar[r]\ar[d] & \pi^*TX \ar[r]\ar@{=}[d] & \pi^*\sO_D(D) \ar[r]\ar[d] & 0 \\
0\ar[r] & TX_{D,r} \ar[r] & \pi^*TX \ar[r] & \pi^*\sO_X(D)\otimes\sO_{\fG_{r-1}} \ar[r] & 0 \\
0 \ar[r] & \sO_{X_{D,r}} \ar[r]\ar[d] & \sT^{\otimes r} \ar@{=}[d]\ar[r] & \pi^*\sO_D(D) \ar[r]\ar[d] & 0 \\
0 \ar[r] & \sT \ar[r] & \sT^{\otimes r} \ar[r] & \pi^*\sO_X(D)\otimes\sO_{\fG_{r-1}} \ar[r] & 0}$$
For the second diagram, note that $\sT^{\otimes r}\cong \pi^*\sO_X(D)$ and that $\fG_{r-1}$ is cut out by a section of $\sT^{\otimes (r-1)}$.  By applying the snake lemma to both diagrams, we see that the cokernel of $\pi^*\sE\to TX_{D,r}$ is isomorphic to the cokernel of $\sO_{X_{D,r}}\to\sT$, which is $\sT\otimes\sO_{\fG}$ since this morphism is the tautological section which cuts out $\fG$.  So we have an exact sequence
\begin{equation}
\label{eq:elmntry_transf}
0\to\pi^*\sE\to TX_{D,r} \to \sT\otimes\sO_{\fG} \to 0,
\end{equation}
and therefore
\begin{equation}
\label{eq:c_1TX}
c_1(TX_{D,r}) = \pi^*(c_1(TX) - \frac{r-1}{r}[D]).
\end{equation}

\begin{prop}
\label{gerbe_exact}
There is an exact sequence of sheaves on $\fG$:
$$0\to\sO_{\fG} \to \pi^*\sE\otimes\sO_{\fG} \to TX_{D,r}\otimes\sO_{\fG} \to \sT\otimes\sO_{\fG} \to 0.$$
\end{prop}

Note that $\pi^*\sE\otimes\sO_{\fG}/\sO_{\fG}$ is locally free since it is the kernel of a surjection of locally free sheaves.
\smallskip

\begin{proof} This follows from (\ref{eq:elmntry_transf}), given the fact that $\Tor_1^{\sO_{X_{D,r}}}(\sT\otimes\sO_{\fG},\sO_{\fG})=\sO_{\fG}$, which follows by tensoring the exact sequence $$0\to \sO_{X_{D,r}}\to \sT\to \sT\otimes\sO_{\fG}\to 0$$ with $\sO_{\fG}$.\end{proof}

\begin{remark}  The above exact sequences all hold if $X$ is a smooth twisted curve and $D$ is an untwisted point of $X$.  Therefore, (\ref{eq:elmntry_transf}) implies that
$$T\fC\cong TC(-\sum p_i)\otimes\prod\sT_i,$$
where $\fC$ is a smooth twisted curve with coarse moduli scheme $C$, $p_i$ are the points of $C$ which are twisted in $\fC$, and $\sT_i$ is the tautological sheaf corresponding to $p_i$.  Recalling the notation from the beginning of section \ref{sec:vfc}.1, it follows that $(\sT_i)_{x_i}$ is the standard representation, where $x_i\in\fC$ is the preimage of $p_i$.
\end{remark}

Let $f:\fC\to X_{D,r}$ be a genus $g$ twisted stable map of class $\beta$ and contact type $\vec{\varrho}$ (Def. \ref{contact_type}).  Then equation \ref{eq:c_1TX} implies that
\begin{equation}
\label{eq:deg_tbdle}
\deg(f^*TX_{D,r})=\deg(f^*\pi^*TX)-\frac{r-1}{r}\beta\cdot D=-\beta\cdot (K_X+D) + \frac{1}{r}\beta\cdot D.
\end{equation}
Proposition \ref{gerbe_exact} implies that the age of $f^*TX_{D,r}$ at a twisted point equals the age of $f^*\sT$ at the same point.  From the definition of contact type, we have
\begin{equation}
\label{eq:age_tbdle}
\age(f^*TX_{D,r},x_i) = \frac{\varrho_i}{r}
\end{equation}
where $x_i$ is the $i$th marked point.  Since $\rk(TX_{D,r})=\rk(TX)=\dim(X)$, Theorem \ref{Riemann-Roch} implies that
\begin{equation}
\label{eq:euchar_tbdle}
\chi(f^*TX_{D,r})=-\beta\cdot (K_X+D) + \frac{1}{r}(\beta\cdot D - \sum_{i=1}^n \varrho_i) - \dim(X)(g-1)
\end{equation}

\subsection{Expected dimension of $\eK_{g,n}(X_{D,r},\beta,\vec{\varrho})$}

Fix $g,n,\beta,\varrho$ and let $\eK=\eK_{g,n}(X_{D,r},\beta,\vec{\varrho})$.  From equation \ref{eq:euchar_tbdle}, it is clear that the restriction of the virtual fundamental class to $\eK$ (which we denote $[\eK]^v$) is homogeneous in the Chow group of $\eK$.  Its degree is the expected dimension of $\eK$ and is denoted $\edim(\eK)$.  This number is important in Gromov-Witten theory because it determines which Gromov-Witten invariants can be nonzero.  For example, suppose we choose for each $i$ an element $a_i\in N^{c_i}((\widehat{X_{D,r}})_{\varrho_i})\subseteq N^*(\widehat{X_{D,r}})$.  Then the Gromov-Witten invariant $I^g_{\beta}(a_1\cdots a_n)$ can be nonzero only if $\edim(\eK)=\sum c_i$.

Since the dimension of $\fM_{g,n}^{\tw}$ is $3g-3+n$, we can compute the expected dimension from equations \ref{eq:edim_gen} and \ref{eq:euchar_tbdle}.
\begin{equation}
\label{eq:edim}
\edim(\eK)=-\beta\cdot (K_X+D)+\frac{1}{r}(\beta\cdot D-\sum_{i=1}^n\varrho_i) + n + (3-\dim(X))(g-1)
\end{equation}

\section{Some invariants of \boldmath $\pd$}
\label{sec:computations}

\subsection{Preliminaries}

In this section we compute some Gromov-Witten invariants of $\pd$, where $D\subseteq\bP^2$ is a smooth curve of degree $\delta$.  In the next section we use these invariants together with associativity of the big quantum product to find recursions that determine all the genus $0$ invariants.

We denote the class of a curve in $N_1(\bP^2)$ by its degree $d$.  Since we are only interested in genus $0$ invariants, we use $I_d$ to mean $I_d^0$.  Recall that $\widehat{\pd}\cong\bP^2\sqcup D$.  We fix the following basis for $N^*(\widehat{\pd})$.

$T_0$ is the unit class of $\bP^2$.

$T_1$ is the hyperplane class of $\bP^2$.

$T_2$ is the class of a point in $\bP^2$.

$T_3$ is the unit class of $D$.

$T_4$ is the class of a point in $D$.
\smallskip

Our first computation is an immediate consequence of Proposition \ref{invts_zero}.  We use $T^{\vec{n}}$ as shorthand for $\prod_{i=0}^4 T_i^{n_i}$.
\begin{align}
I_0(T^{\vec{n}})=0 & \mbox{ if } & n_2+n_4\ge 2 \label{eq:deg_0_zero} \\
I_1(T^{\vec{n}})=0 & \mbox{ if } & n_2+n_4\ge 3 \label{eq:deg_1_zero}
\end{align}

\subsection{Degree 0, 3-point invariants}

For any triple $\vec{\varrho}$ of integers $0$ and $1$, let $\eK_{\vec{\varrho}}=\eK_{0,3}(\pd,0,\vec{\varrho})$.  By equation \ref{eq:edim}, the expected dimension of $\eK_{\vec{\varrho}}$ is
\begin{equation}
\label{eq:edimK_rho}
2-\frac{\varrho_1+\varrho_2+\varrho_3}{2}.
\end{equation}
Since the expected dimension must be an integer, it follows from this that $\eK_{\vec{\varrho}}=\emptyset$ if there are an odd number of $1$'s in $\vec{\varrho}$.  If $\vec{\sigma}$ is any permutation of $\vec{\varrho}$, then clearly $\eK_{\vec{\varrho}}\cong\eK_{\vec{\sigma}}$, so there are essentially two distinct cases.  Note that there is a natural morphism $F_{\vec{\varrho}}:\eK_{\vec{\varrho}}\to\overline{\eM}_{0,3}(\bP^2,0)\cong\bP^2$ which equals the composition of any evaluation map $\eK_{\vec{\varrho}}\to \widehat{\pd}$ with the projection $\widehat{\pd}\to\bP^2$.

It follows from \cite[4.2]{stacks} that $F_{(0,0,0)}$ is an isomorphism over $\bP^2\setminus D$.  Since the expected dimension of $\eK_{(0,0,0)}$ is 2 and $\bP^2\setminus D$ is smooth of dimension 2, it follows from Proposition \ref{vfc_coeff=1} that $(F_{(0,0,0)})_*[\eK_{(0,0,0)}]^v=[\bP^2]$.  It follows that
\begin{equation}
I_0(T_0^2T_2)=I_0(T_0T_1^2)=1 \label{eq:deg_0_first}
\end{equation}
and that all other degree $0$, $3$-point invariants involving neither $T_3$ nor $T_4$ are zero.

Now we study $\eK_{(1,1,0)}$.  Since $e_1$ and $e_2$ map to $D$, $F_{(1,1,0)}$ factors through $D$.  For any closed point $x\in D$, we now describe the unique point in the fiber of $F_{(1,1,0)}$ over $x$.  Let $\fC$ be the square root of $\bP^1$ at two distinct points, let $\sT_1$ and $\sT_2$ be the tautological sheaves corresponding to the resulting twisted points, and let $\sT$ be the tautological sheaf on $\pd$.  There is a morphism $f:\fC\to\pd$ which on the level of coarse moduli spaces sends $\bP^1$ to $x$ and satisfies $f^*\sT\cong\sT_1\otimes\sT_2\otimes\sO_{\fC}(-1)$.  Since $\sT_1$ and $\sT_2$ appear with exponent 1, this morphism has contact type $(1,1,0)$.  It is easy to see that any stable map in the fiber over $x$ has such an underlying morphism $f$, and also that there is a unique such stable map up to isomorphism.  Moreover, $f:\fC\to\pd$ has a unique nontrivial 2-automorphism given by multiplication by $-1$ on $f^*\sT$.

We can now compute the deformation and obstruction spaces for this stable map.  By Proposition \ref{gerbe_exact}, we have an exact sequence
$$0 \to \sO_{\fC} \to f^*TX_{D,r} \to \sT_1\otimes\sT_2\otimes\sO_{\fC}(-1) \to 0.$$
If $\gamma:\fC\to\bP^1$ is the projection, then it follows from \cite[3.1.1]{stacks} that $\gamma_*$ of the right hand term is $\sO_{\bP^1}(-1)$, so it has no cohomology.  Hence $h^0(\fC,f^*TX_{D,r})=1$ and $h^1(\fC,f^*TX_{D,r})=0$.  Since this holds for every stable map in $\eK_{(1,1,0)}$, it follows from \cite[7.3]{BF} that $\eK_{(1,1,0)}$ is smooth of the expected dimension and $[\eK_{(1,1,0)}]^v=[\eK_{(1,1,0)}]$.

Let $U\to\eK_{(1,1,0)}$ be an \'etale surjective map.  Then $U$ is smooth and 1-dimensional and every component of $U$ dominates $D$.  It follows by \cite[III-9.7]{Ha} that $U\to D$ is flat, which by definition implies that $F_{(1,1,0)}$ is flat.  Since $\eK_{(1,1,0)}$ is proper, $F_{(1,1,0)}$ is also proper and hence finite.  The degree of $F_{(1,1,0)}$ is the degree of $F_{(1,1,0)}^{-1}(x)$ for any $x\in D$.  By generic smoothness, the preimage of a general $x\in D$ is reduced, so it follows by the above analysis that $F_{(1,1,0)}^{-1}(x)\cong B\mu_2$ for general $x\in D$, and hence $F_{(1,1,0)}$ has degree $1/2$.  Therefore, $(F_{(1,1,0)})_*[\eK_{(1,1,0)}]^v= (1/2)[D]$.

It now follows that the only nonzero degree $0$, $3$-point invariants involving either $T_3$ or $T_4$ are
\begin{equation}
\label{eq:deg_0_second}
I_0(T_1T_3^2) = \delta/2 \mbox{ and } I_0(T_0T_3T_4)=1/2.
\end{equation}

\subsection{Some degree 1 invariants}

Now we compute the invariants
\begin{equation}
\label{eq:deg_1_invts}
I_1(T_2^2 T_3^{\delta})=\delta ! \mbox{ and } I_1(T_2 T_3^{\delta-1} T_4)=(\delta-1)!.
\end{equation}
It follows in the same way that $I_1(T_3^{\delta-2} T_4^2)=(\delta-2)!$ if $\delta\ge 2$, but this invariant also follows from the associativity relations.

Let $\varrho_i=0$ for $i=1,2$ and $\varrho_i=1$ for $3\le i\le \delta + 2$.  We have a finite morphism $F:\eK_{0,\delta+2}(\pd,1,\vec{\varrho})\to\overline{\eM}_{0,\delta+2}(\bP^2,1)$ by \cite[1.4.1]{AV}.  Since this is compatible with the evaluation maps, it factors through $\cap_{i=3}^{\delta+2}e_i^{-1}(D)\subseteq\overline{\eM}_{0,\delta+2}(\bP^2,1)$.

Let $p:\overline{\eM}_{0,\delta+2}(\bP^2,1)\to\overline{\eM}_{0,2}(\bP^2,1)$ be the flat and proper morphism which forgets the last $\delta$ markings.  Let $U\subseteq\overline{\eM}_{0,2}(\bP^2,1)$ be the dense open subscheme consisting of stable maps $f:\bP^1\to\bP^2$ such that $f(\bP^1)$ is transverse to $D$ and the marked points do not map into $D$.  Let $G$ and $H$ be as in the diagram below.
$$\xymatrix{
\eK_{0,\delta+2}(\pd,1,\vec{\varrho}) \ar[r] \ar[dr]^{G} \ar@/^1pc/[rr]^{F} & \cap_{i=3}^{\delta +2}e_i^{-1}(D) \ar@{^{(}->}[r] \ar[d]^{H} & \overline{\eM}_{0,\delta +2}(\bP^2,1) \ar[dl]_{p} \\
U \ar@{^{(}->}[r] & \overline{\eM}_{0,2}(\bP^2,1)}$$

If $f:\fC\to\pd$ is a stable map in $G^{-1}(U)$, then we claim that $\fC$ is smooth.  If $\fC$ had a node, then there would be a component mapping with degree $0$ which would have to contain at least two marked points.  Both would have to be twisted by definition of $U$.  If $\fC_0\subseteq\fC$ is the irreducible component which maps with positive degree, then the contact type of $f\vert_{\fC_0}$ must be odd at every preimage of $\fG$ by \cite[3.3.6]{stacks}.  This implies that every such point must be twisted, but since there are $\delta$ points in the preimage of $\fG$ and only $\delta$ twisted markings, there cannot be a node.

Let $V\subseteq\cap_{i=3}^{\delta+2} e_i^{-1}(D)$ be the open subscheme containing only maps from smooth curves into $\bP^2$.

\begin{prop}
\label{isom_over_U}
$G^{-1}(U)\to H^{-1}(U)\cap V$ is an isomorphism.
\end{prop}

\begin{proof} This follows from \cite[4.2]{stacks} after observing that, in the notation used there, $\eU_{0,\delta + 2}(\pd,1,\vec{\varrho})=G^{-1}(U)$ and $\eV_{0,\delta + 2}(\bP^2,1,\vec{\varrho})=H^{-1}(U)\cap V$.\end{proof}

\begin{prop}
\label{finite_etale}
$H^{-1}(U)\cap V\to U$ is finite and \'etale of degree $\delta !$.
\end{prop}

\begin{proof} To show the morphism is \'etale, we use the following criterion.  Let $R$ be a Noetherian ring and $I\subseteq R$ a nilpotent ideal.  Given a commutative diagram
$$\xymatrix{
\spec (R/I) \ar[r] \ar[d] & H^{-1}(U)\cap V \ar[d] \\
\spec R \ar[r] & U,}$$
it must be shown that there is a unique morphism $\spec R\to H^{-1}(U)\cap V$ making the diagram commute.  This means we have a smooth family of rational curves $C\to\spec R$ with two disjoint sections $s_1,s_2$ and a morphism $f:C\to\bP^2$ such that over $\spec (R/I)$, $f^{-1}(D)$ is a disjoint union of $\delta$ sections which are disjoint from $s_1$ and $s_2$.  Since $I$ is nilpotent, $f^{-1}(D)$ must have $\delta$ connected components.  Each is clearly finite over $\spec R$, and each is flat by \cite[5.2]{stacks}.  Since a flat, finite, degree 1 morphism is an isomorphism, each connected component of $f^{-1}(D)$ is a section.  This verifies the criterion.

The morphism is proper by Proposition \ref{isom_over_U} since $G$ is proper.  To show it is finite of degree $\delta !$, it suffices to show that each closed point has $\delta !$ preimages.  The closed point is represented by a line transverse to $D$ with two markings away from $D$, and a preimage is determined by an ordering of the $\delta$ intersection points.  This completes the proof.\end{proof}

Note that the first two evaluation maps $\eK:=\eK_{0,\delta+2}(\pd,1,\vec{\varrho})\to\bP^2$ factor through $G$.  Since the product of evaluation maps $\overline{\eM}_{0,2}(\bP^2,1)\to\bP^2\times\bP^2$ is a birational morphism, it follows from the next proposition that $I_1(T_2^2 T_3^{\delta}) = \delta!$.

\begin{prop}
\label{virt_pushforward}
$G_*[\eK]^v=\delta ! [\overline{\eM}_{0,2}(\bP^2,1)]$.
\end{prop}

\begin{proof} By equation \ref{eq:edim}, $[\eK]^v$ is homogeneous of degree $4$ in the Chow group of $\eK$.  Since $U$ is dense in $\overline{\eM}_{0,2}(\bP^2,1)$, the only contribution to $G_*[\eK]^v$ comes from the closure of $G^{-1}(U)$ in $\eK$.  By Propositions \ref{isom_over_U} and \ref{finite_etale}, $G^{-1}(U)$ is smooth of the expected dimension, so the result now follows from Propositions \ref{vfc_coeff=1} and \ref{finite_etale}.\end{proof}

To compute $I_1(T_2 T_3^{\delta-1} T_4)$, we work instead with the following diagram.
$$\xymatrix{
\eK_{0,\delta +1}(\pd,1,\vec{\sigma}) \ar[r] \ar[dr]^G & \cap_{i=2}^{\delta + 1} e_i^{-1}(D) \ar@{^{(}->}[r] \ar[d]^H & \overline{\eM}_{0,\delta + 1}(\bP^2,1) \ar[d] \\
U \ar@{^{(}->}[r] & e_2^{-1}(D) \ar@{^{(}->}[r] & \overline{\eM}_{0,2}(\bP^2,1)}$$
Here $\sigma_1=0$ and $\sigma_i=1$ for $2\le i\le\delta +1$.  The dense open set $U\subseteq e_2^{-1}(D)$ consists of lines transverse to $D$ where the first marked point lies off of $D$.  Let $V\subseteq \cap_{i=2}^{\delta+1} e_i^{-1}(D)$ be the open subscheme of maps from smooth curves into $\bP^2$.

By the same arguments as above, one can show that $G^{-1}(U)\to H^{-1}(U)\cap V$ is an isomorphism, that $H^{-1}(U)\cap V\to U$ is finite and \'etale of degree $(\delta-1)!$, and that the morphism $U\to\overline{\eM}_{0,1}(\bP^2,1)$ forgetting the second marking is \'etale.  It follows that $G^{-1}(U)$ is smooth of the expected dimension and that $G_*[\eK]^v=(\delta - 1)![e_2^{-1}(D)]$.  Then the result follows from the fact that the exceptional locus of $\overline{\eM}_{0,2}(\bP^2,1)\to\bP^2\times\bP^2$ does not contain $\bP^2\times D$.

\section{The big quantum product}
\label{sec:quant_prod}

\subsection{Preliminaries}

We continue with the notation of the last section.  So $D\subseteq\bP^2$ is a smooth curve of degree $\delta$ and $T_0,\ldots,T_4$ is the chosen basis of $N^*(\widehat{\pd})$.

First we define the big quantum product for $\pd$ and then we use the fact that it is associative to compute recursions.  While the big quantum product exists for any smooth Deligne-Mumford stack having projective coarse moduli scheme, we only define it for $\pd$ in order to simplify the notation.  The reader may also wish to consult \cite{AGV,CR}.  We have adopted the notation of \cite[\S 8]{FP}.

Let $(g^{ij})_{0\le i,j\le 4}$ be the matrix
\def\arraystretch{1}
$$\left(\begin{array}{ccccc}
0 & 0 & 1 & 0 & 0\\
0 & 1 & 0 & 0 & 0\\
1 & 0 & 0 & 0 & 0\\
0 & 0 & 0 & 0 & 2\\
0 & 0 & 0 & 2 & 0
\end{array}\right).$$
\def\arraystretch{1.6}
This corresponds to the inverse matrix of the intersection form
$$a\otimes b\mapsto \int_{\fI(\pd)}a\cdot b$$
expressed in the basis $\{\pi^*T_i\}$, where $\fI(\pd)$ is the inertia stack of $\pd$ and $\pi:\fI(\pd)\to\widehat{\pd}$ is the projection.  The reason for the 2's is that one gets a factor of $1/2$ when one integrates over $\fG$ instead of $D$.  In general, an involution on the inertia stack enters into the definition of $g^{ij}$ but for $\pd$ this involution is the identity.

The definition of the big quantum product uses a generating function for the genus $0$ Gromov-Witten invariants $I_d(T_0^{n_0}\cdots T_4^{n_4})$ called the quantum potential.  It is a power series in indeterminates $y_0,\ldots,y_4$ given by
$$\Phi(y_0,\ldots,y_4)=\sum_{n_0+\cdots +n_4\ge 3}\sum_{d=0}^\infty q^d I_d(T_0^{n_0}\cdots T_4^{n_4})\frac{y_0^{n_0}}{n_0!}\cdots\frac{y_4^{n_4}}{n_4!}.$$ Let
$$\Phi_{ijk}=\frac{\partial^3\Phi}{\partial y_i\partial y_j\partial y_k}$$
for $0\le i,j,k\le 4.$

\begin{definition}
\label{quantum_product}
The big quantum product is the $R:=\bQ[[y_0,\ldots,y_4,q]]$-linear product on the free $R$-module with basis $T_0,\ldots,T_4$ which is given by
$$T_i * T_j = \sum_{e,f=0}^4 \Phi_{ije}g^{ef}T_f.$$
\end{definition}

\subsection{Forgetting an untwisted point}

Here we apply the forgetting a point axiom \cite{AGV2} to the stacks $\bP^2_{D,r}$.  Let $\vec{\varrho}$ be an $n$-tuple and let $\vec{\sigma}=(\varrho_1,\ldots,\varrho_n,0)$.  If either $n\ge 3$ or $d>0$, then it follows from \cite[9.1.3]{AV} that there is a morphism $F:\eK_{0,n+1}(\pd,d,\vec{\sigma})\to\eK_{0,n}(\pd,d,\vec{\varrho})$ which forgets the last marked point.  The forgetting a point axiom says firstly that $\eK_{0,n+1}(\pd,d,\vec{\sigma})$ is isomorphic to the universal curve over $\eK_{0,n}(\pd,d,\vec{\varrho})$ in such a way that $F$ is the projection and $e_{n+1}$ is the composition of the univeral morphism with the projection $\pd\to\bP^2$.  The second part of the axiom is that $$F^*[\eK_{0,n}(\pd,d,\vec{\varrho})]^v=[\eK_{0,n+1}(\pd,d,\vec{\sigma})]^v.$$

From this, the following equations can be derived almost exactly as in the case of ordinary stable maps (cf. \cite[\S 7, I-III]{FP}).

\begin{enumerate}
   \item If $n_0>0$ and either $d>0$ or $\sum n_i>3$, then
      \begin{equation}\label{eq:GW_relations1}
            I_d(T^{\vec{n}})=0.\end{equation}
   \item If $n_1>0$ and either $d>0$ or $\sum n_i>3$, then
      \begin{equation}\label{eq:GW_relations2}
            I_d(T^{\vec{n}})=dI_d(T_0^{n_0}T_1^{n_1-1}T_2^{n_2}T_3^{n_3}T_4^{n_4}).\end{equation}
   \item If $n_2>0$ and $\sum n_i>3$, then
      \begin{equation}\label{eq:GW_relations3}
            I_0(T^{\vec{n}})=0.\end{equation}
\end{enumerate}

\subsection{Identity and associativity}

The quantum product is clearly commutative.  We also need the fact that $T_0$ is a multiplicative identity and that the product is associative.

\begin{thm}
\label{identity}
For $0\le i\le 4$, $T_0*T_i=T_i$.
\end{thm}

\begin{proof} By equation \ref{eq:GW_relations1}, the only nonzero invariants with $n_0>0$ are those of the form $I_0(T_0T_iT_e)$ for any $i$ and $e$.  From definition \ref{quantum_product}, we see that $T_0*T_i=\sum I_0(T_0T_iT_e)g^{ef}T_f$.  The invariants $I_0(T_0T_iT_e)$ are computed in equations \ref{eq:deg_0_first} and \ref{eq:deg_0_second}.  The theorem now follows from the definition of $g^{ij}$.\end{proof}

For associativity of the quantum product, see \cite{AGV2}:

\begin{equation}\label{eq:associativity}
(T_i*T_j)*T_k=T_i*(T_j*T_k) \mbox{ for all } i,j,k.
\end{equation}

\subsection{Simplifications}

Let $\Psi$, $\Psi'$, and $\Gamma$ be the power series defined by the following formulas.

\begin{align}
\Psi & = \frac{1}{6}\sum_{i,j,k=0}^4 I_0(T_iT_jT_k)y_iy_jy_k \notag \\
\Psi' & = \sum_{n_3+n_4\ge 4} I_0(T_3^{n_3}T_4^{n_4})\frac{y_3^{n_3}}{n_3!}\frac{y_4^{n_4}}{n_4!} \label{eq:quant_pot_comps} \\
\Gamma & = \sum_{n_2+n_3+n_4\ge 0}\sum_{d=1}^\infty (qe^{y_1})^d I_d(T_2^{n_2}T_3^{n_3}T_4^{n_4})\frac{y_2^{n_2}}{n_2!}\frac{y_3^{n_3}}{n_3!}\frac{y_4^{n_4}}{n_4!} \notag
\end{align}

By equations \ref{eq:GW_relations1}-\ref{eq:GW_relations3}, $\Psi+\Psi'+\Gamma$ is congruent to $\Phi$ modulo terms of degree less than or equal to $2$.  Since the big quantum product only involves the third order partial derivatives of $\Phi$, we can use $\Psi+\Psi'+\Gamma$ in place of $\Phi$ in Definition \ref{quantum_product}.

We now introduce the stringy product, denoted $\cdot_s$.  For our purposes, it is just a convenient way to encode the effect of $\Psi$ on the big quantum product.  It is defined by $$T_i\cdot_s T_j=\sum_{e,f}\Psi_{ije}g^{ef}T_f.$$  From equation \ref{eq:deg_0_second}, we compute $T_3\cdot_s T_3=(\delta/2)T_1$, $T_3\cdot_s T_1=\delta T_4$, $T_3\cdot_s T_4=(1/2)T_2$, and $T_3\cdot_s T_2=0$.  By associativity, commutativity, and the fact that $T_0$ is a unit, this determines the stringy product.

Once we have the stringy product, we only need to use $\Psi'$ and $\Gamma$, and these only involve $T_2$, $T_3$, and $T_4$.  So we use $I_d(n_2,n_3,n_4)$, or sometimes $I_d(\vec{n})$, to denote $I_d(T_2^{n_2}T_3^{n_3}T_4^{n_4})$.  It is important to know which invariants can be nonzero.

\begin{prop}
\label{nonzero_invts}
If $I_d(n_2,n_3,n_4)\ne 0$, then $$3d-1=\frac{d\delta + n_4-n_3}{2} + n_2.$$
\end{prop}

\begin{proof} Let $\varrho_i=0$ for $1\le i\le n_2$ and $\varrho_i=1$ for $n_2+1\le i\le n_2+n_3+n_4$.  Equation \ref{eq:edim} implies that the expected dimension of $\eK_{0,n_2+n_3+n_4}(\pd, d, \vec{\varrho})$ is $3d-d\delta/2 + (n_3+n_4)/2 + n_2 - 1$.  The equality comes from setting this equal to $2n_2+n_4$.\end{proof} 

If we apply this to $\Psi'$, we see that $I_0(0,n_3,n_4)$ can only be nonzero when $n_3-n_4=2$.  Equation \ref{eq:deg_0_zero} imposes the additional condition $n_4\le 1$.  We have thus shown that $\Psi'=\lambda y_3^3y_4/6$, where $\lambda:=I_0(0,3,1)$.

It follows that the quantum product is given by
\begin{align}
T_i*T_j & = T_i\cdot_s T_j + \sum_{ef} (\Psi'_{ije}+\Gamma_{ije})g^{ef}T_f \notag \\
 & = T_i\cdot_s T_j + \Gamma_{ij1}T_1 + \Gamma_{ij2}T_0 + 2\Gamma_{ij3}T_4 + 2\Gamma_{ij4}T_3 + \lambda \frac{\partial^2 y_3^2 y_4}{\partial y_i\partial y_j} T_4 + \frac{\lambda}{3}\frac{\partial^2 y_3^3}{\partial y_i\partial y_j} T_3. \label{eq:big_quant_prod_simpl}
\end{align}

\subsection{Recursions}

Using equation \ref{eq:big_quant_prod_simpl} together with associativity, it is a simple but tedious computation to compute the recursions.  The only difficulty is in knowing which products to apply associativity to, and we do not claim to have found the most efficient algorithm.  The following four relations are obtained by comparing respectively the coefficients of $T_0$ in $(T_1*T_1)*T_2$ and $T_1*(T_1*T_2)$, those of $T_3$ in $(T_3*T_3)*T_4$ and $T_3*(T_3*T_4)$, those of $T_3$ in $(T_3*T_1)*T_4$ and $T_3*(T_1*T_4)$, and those of $T_1$ in $(T_3*T_3)*T_1$ and $T_3*(T_3*T_1)$.

\begin{gather}
\Gamma_{222} = \Gamma_{112}^2-\Gamma_{111}\Gamma_{122}+2(2\Gamma_{123}\Gamma_{124}-\Gamma_{113}\Gamma_{224}-\Gamma_{114}\Gamma_{223}) \label{relation1} \\
\delta\Gamma_{144}+4\lambda (y_4\Gamma_{444}-y_3\Gamma_{344})-2\Gamma_{234} = 2(\Gamma_{134}^2-\Gamma_{133}\Gamma_{144})+4(\Gamma_{334}\Gamma_{344}-\Gamma_{333}\Gamma_{444}) \label{relation2} \\
2\delta\Gamma_{444}-\Gamma_{124}-4\lambda y_3\Gamma_{144} = 2(\Gamma_{114}\Gamma_{134}-\Gamma_{113}\Gamma_{144})+4(\Gamma_{144}\Gamma_{334}-\Gamma_{133}\Gamma_{444}) \label{relation3} \\
\begin{split}
\Gamma_{233}+\delta(\frac{1}{2}\Gamma_{111}-2\Gamma_{134})+2\lambda & (y_3\Gamma_{113}+ y_4\Gamma_{114}) = \\
& \Gamma_{113}^2-\Gamma_{111}\Gamma_{133}+2(2\Gamma_{133}\Gamma_{134}-\Gamma_{113}\Gamma_{334}-\Gamma_{114}\Gamma_{333}) 
\end{split}
\label{relation4}
\end{gather}

It is now possible to compute $\lambda$ using equation \ref{relation4}.  Comparing the coefficients of $qe^{y_1}y_2^2y_3^{\delta}/2\delta !$ yields
$$\delta(\frac{1}{2}+2\lambda)I_1(2,\delta,0)=2\delta I_1(2,\delta+1,1)-I_1(3,\delta+2,0).$$
The right hand side is zero by equation \ref{eq:deg_1_zero} and $I_1(2,\delta,0)$ is nonzero by equation \ref{eq:deg_1_invts}.  Therefore, $\lambda=-1/4$.

By comparing coefficients in equations \ref{relation1}-\ref{relation4}, we obtain recursions \ref{recursion1}-\ref{recursion4} in that order.  In each recursion, $d_1$ and $d_2$ vary over positive integers and $\vec{p}:=(p_2,p_3,p_4)$ and $\vec{q}:=(q_2,q_3,q_4)$ vary over triples of nonnegative integers.  In addition to the condition in parentheses, we also assume $d>0$.

\begin{multline}
 \label{recursion1} 
I_d(n_2,n_3,n_4)= \\
(n_2\ge 3) \hfill \sum_{\stackrel{\scriptstyle d_1+d_2=d}{\vec{p}+\vec{q}=\vec{n}-(1,0,0)}}I_{d_1}(\vec{p})I_{d_2}(\vec{q})\ch{n_3}{p_3}\ch{n_4}{p_4}\left[d_1^2d_2^2\ch{n_2-3}{p_2-1}-d_1^3d_2\ch{n_2-3}{p_2}\right]+ \\ 
\sum_{\stackrel{\scriptstyle d_1+d_2=d}{\vec{p}+\vec{q}=\vec{n}+(-1,1,1)}}2I_{d_1}(\vec{p})I_{d_2}(\vec{q})\ch{n_3}{p_3-1}\ch{n_4}{p_4}\left[2d_1d_2\ch{n_2-3}{p_2-1}-d_1^2\ch{n_2-3}{p_2}-d_2^2\ch{n_2-3}{p_2-2}\right]
\end{multline}
\begin{multline}
\label{recursion2}
(d\delta+n_3-n_4+2)I_d(n_2,n_3,n_4)=2I_d(n_2+1,n_3+1,n_4-1)+ \\
(n_4\ge 2) \hfill \sum_{\stackrel{\scriptstyle d_1+d_2=d}{\vec{p}+\vec{q}=\vec{n}+(0,2,0)}}2d_1d_2I_{d_1}(\vec{p})I_{d_2}(\vec{q})\ch{n_2}{p_2}\left[\ch{n_3}{p_3-1}\ch{n_4-2}{p_4-1}-\ch{n_3}{p_3-2}\ch{n_4-2}{p_4}\right]+ \\
\sum_{\stackrel{\scriptstyle d_1+d_2=d}{\vec{p}+\vec{q}=\vec{n}+(0,3,1)}}4I_{d_1}(\vec{p})I_{d_2}(\vec{q})\ch{n_2}{p_2}\left[\ch{n_3}{p_3-2}\ch{n_4-2}{p_4-1}-\ch{n_3}{p_3-3}\ch{n_4-2}{p_4}\right]
\end{multline}
\begin{multline}
\label{recursion3}
2\delta I_d(n_2,n_3,n_4)=dI_d(n_2+1,n_3,n_4-2)-n_3dI_d(n_2,n_3-1,n_4-1)+ \\
(n_4\ge 3) \hfill \sum_{\stackrel{\scriptstyle d_1+d_2=d}{\vec{p}+\vec{q}=\vec{n}+(0,1,-1)}}2I_{d_1}(\vec{p})I_{d_2}(\vec{q})\ch{n_2}{p_2}\ch{n_3}{p_3-1}\left[d_1d_2^2\ch{n_4-3}{p_4-1}-d_1^2d_2\ch{n_4-3}{p_4}\right]+ \\
\sum_{\stackrel{\scriptstyle d_1+d_2=d}{\vec{p}+\vec{q}=\vec{n}+(0,2,0)}}4I_{d_1}(\vec{p})I_{d_2}(\vec{q})\ch{n_2}{p_2}\ch{n_3}{p_3-2}\left[d_2\ch{n_4-3}{p_4-1}-d_1\ch{n_4-3}{p_4}\right]
\end{multline}
\begin{align}
\frac{1}{2}d^2(d\delta-n_3-n_4)I_d(n_2,n_3,n_4)=2d\delta I_d(n_2,n_3+1,n_4+1)-I_d(n_2+1,n_3+2,n_4) & + \notag \\
\sum_{\stackrel{\scriptstyle d_1+d_2=d}{\vec{p}+\vec{q}=\vec{n}+(0,2,0)}}I_{d_1}(\vec{p})I_{d_2}(\vec{q})\ch{n_2}{p_2}\ch{n_4}{p_4}\left[d_1^2d_2^2\ch{n_3}{p_3-1}-d_1^3d_2\ch{n_3}{p_3}\right] & +  \label{recursion4} \\
\sum_{\stackrel{\scriptstyle d_1+d_2=d}{\vec{p}+\vec{q}=\vec{n}+(0,3,1)}}2I_{d_1}(\vec{p})I_{d_2}(\vec{q})\ch{n_2}{p_2}\ch{n_4}{p_4}\left[2d_1d_2\ch{n_3}{p_3-2}-d_1^2\ch{n_3}{p_3-1}-d_2^2\ch{n_3}{p_3-3}\right] & \notag
\end{align}

Combining these recursions with equation \ref{eq:deg_1_invts} and Proposition \ref{nonzero_invts}, the genus $0$ Gromov-Witten invariants of $\pd$ can be computed by the following algorithm.

\begin{algorithm}
Let $d\ge 1$ and $n_2,n_3,n_4\ge 0$ be integers.  Compute $I:=I_d(n_2,n_3,n_4)$ as follows.
\begin{enumerate}
  \item If $3d-1\ne (d\delta+n_4-n_3)/2+n_2$, then $I=0$.
  \item Otherwise, if $(d,n_2,n_3,n_4)=(1,2,\delta,0)$, then $I=\delta !$.
  \item Otherwise, if $(d,n_2,n_3,n_4)=(1,1,\delta-1,1)$, then $I=(\delta-1)!$.
  \item Otherwise, if $n_2\ge 3$, apply recursion \ref{recursion1}.
  \item Otherwise, if $n_4\ge 2$ and $n_4-n_3\ne d\delta + 2$, apply recursion \ref{recursion2}.
  \item Otherwise, if $n_4-n_3=d\delta + 2$, apply recursion \ref{recursion3}.
  \item Otherwise, apply recursion \ref{recursion4}.
\end{enumerate}
\end{algorithm}

\begin{proof}[Justification]
To justify the algorithm, we first show that if $n_4-n_3=d\delta + 2$ then $n_4\ge 3$, and if step 7 is reached then $n_3+n_4\ne d\delta$.  The first statement follows from $d\delta\ge 1$ and $n_3\ge 0$.  To reach the last step, we must have $n_2\le 2$, $n_4\le 1$, and $3d-1=(d\delta+n_4-n_3)/2 + n_2$.  If $n_3+n_4=d\delta$, then $n_3\ge d\delta-1$, so $3d-1\le 3$.  This implies $d=1$, and the only two possibilities for $(n_2,n_3,n_4)$ are handled in steps 2 and 3.

Since the recursions don't clearly reduce degree $d$ invariants to degree $d-1$ invariants, it must also be shown that the algorithm terminates.  To do this, fix a degree $d$ and define a sequence of triples $(n_2^{(i)}, n_3^{(i)}, n_4^{(i)})$ to be \emph{admissible} if in order to compute $I_d(n_2^{(i)},n_3^{(i)},n_4^{(i)})$, the algorithm requires one to compute $I_d(n_2^{(i+1)}, n_3^{(i+1)}, n_4^{(i+1)})$.  It suffices to show that any admissible sequence starting at a given point has bounded length.

One can see this geometrically by working in the $(n_3,n_4)$ plane, noting that $n_2$ is determined by $n_3$ and $n_4$.  From a given point, one is only allowed to move by five vectors:  $(1,-1)$, $(0,-2)$, $(-1,-1)$, $(1,1)$, and $(2,0)$.  So $n_4-n_3$ is nonincreasing in an admissible sequence.  Moreover, a move by $(-1,-1)$ is only allowed on the line $n_4-n_3=d\delta + 2$, a move by $(1,1)$ is only allowed in the range $n_4\le 1$, and these two sets are disjoint.  Otherwise, $n_4-n_3$ decreases by $2$.

Once an admissible path reaches the line $n_4-n_3=(6-\delta)d-8$, it must terminate since then $n_2=3$.  From this it is easy to get a bound on the length of an admissible path starting at $(n_3,n_4)$.  For example $\frac{1}{2}[n_4-n_3+8-(6-\delta)d]+n_4+\mathrm{max}(0,10-(6-\delta)d)$ works.
\end{proof}

\noindent University of Michigan \\
2074 East Hall \\
Ann Arbor, MI 48109-1043 \\
\ttfamily{cdcadman@umich.edu}

\begin{thebibliography}{XXX}
{
\bibitem[ACV]{ACV}
{\sc D. Abramovich, A. Corti, and A. Vistoli},
Twisted bundles and admissible covers,
{\sl Comm.\ Algebra }31 (2003) 3547--3618.

\bibitem[AGV1]{AGV}
{\sc D. Abramovich, T. Graber, and A. Vistoli},
Algebraic orbifold quantum products,
{\sl Orbifolds in mathematics and physics (Madison, WI, 2001)}, 1--24,
Amer.\ Math.\ Soc., 2002.

\bibitem[AGV2]{AGV2}
{\sc D. Abramovich, T. Graber, and A. Vistoli},
Gromov-Witten theory of Deligne-Mumford stacks,
math.AG/0603151.

\bibitem[AV]{AV}
{\sc D. Abramovich and A. Vistoli},
Compactifying the space of stable maps,
{\sl J.\ Amer.\ Math.\ Soc.\ }15 (2002) 27--75.

\bibitem[Be]{Be}
{\sc K. Behrend},
Gromov-Witten invariants in algebraic geometry,
{\sl Invent.\ Math.\ }127 (1997) 601--617.

\bibitem[BF]{BF}
{\sc K. Behrend and B. Fantechi},
The intrinsic normal cone,
{\sl Invent.\ Math.\ }128 (1997) 45--88.

\bibitem[Ca1]{stacks}
{\sc C. Cadman},
Using stacks to impose tangency conditions on curves, 
{\sl Amer.\ J.\ Math.\ }(to appear), http://www.math.lsa.umich.edu/\verb{~{cdcadman/research/stacks.pdf .

\bibitem[Ca2]{enumerate}
{\sc C. Cadman},
On the enumeration of rational plane curves with tangency conditions,
arXiv: math.AG/0509671.

\bibitem[CC]{enum2}
{\sc C. Cadman and L. Chen},
Enumeration of rational plane curves tangent to a smooth cubic,
in preparation.

\bibitem[CH]{CH}
{\sc L. Caporaso and J. Harris},
Counting plane curves of any genus,
{\sl Invent.\ Math.\ }131 (1998) 345--392.

\bibitem[CR]{CR}
{\sc W. Chen and Y. Ruan},
Orbifold Gromov-Witten theory,
{\sl Orbifolds in mathematics and physics (Madison, WI, 2001)}, 25--85,
Amer.\ Math.\ Soc., 2002.

\bibitem[FP]{FP}
{\sc W. Fulton and R. Pandharipande},
Notes on stable maps and quantum cohomology,
{\sl Algebraic geometry (Santa Cruz, 1995)}, 45--96,
Amer.\ Math.\ Soc., 1997.

\bibitem[GP]{GP}
{\sc T. Graber and R. Pandharipande},
Localization of virtual classes,
{\sl Invent.\ Math.\ }135 (1999) 487--518.

\bibitem[Ha1]{Ha}
{\sc R. Hartshorne},
{\sl Algebraic Geometry},
Springer-Verlag, 1977.

\bibitem[Ha2]{Ha2}
{\sc R. Hartshorne},
Lectures on deformation theory, 
lecture notes, 2004.

\bibitem[Il]{Il}
{\sc L. Illusie},
{\sl Complexe cotangent et d\'eformations I, II},
Springer-Verlag, 1971.

\bibitem[KM]{KM}
{\sc M. Kontsevich and Y. Manin},
Gromov-Witten classes, quantum cohomology, and enumerative geometry,
{\sl Comm.\ Math.\ Phys.\ }164 (1994) 525--562.

\bibitem[LMB]{LMB}
{\sc G. Laumon and L. Moret-Bailly},
{\sl Champs alg\'{e}briques},
Springer-Verlag, 2000.

\bibitem[LT]{LT}
{\sc J. Li and G. Tian},
Virtual moduli cycles and Gromov-Witten invariants of algebraic varieties,
{\sl J.\ Amer.\ Math.\ Soc.\ }11 (1998) 119--174.

\bibitem[Ol]{Ol}
{\sc M. Olsson},
On (log) twisted curves,
preprint.

\bibitem[OS]{OS}
{\sc M. Olsson and J. Starr},
Quot functors for Deligne-Mumford stacks,
{\sl Comm.\ Algebra\ }31 (2003) 4069--4096.

\bibitem[Vi]{Vi}
{\sc A. Vistoli},
Intersection theory on algebraic stacks and on their moduli spaces,
{\sl Invent.\ Math.\ }97 (1989) 613--670.
}
\end{thebibliography}
\end{document}